\documentclass{ifacconf}

\usepackage{graphicx}      
\usepackage{natbib}        
\usepackage{amssymb}
\usepackage{float}
\usepackage{subcaption}
\usepackage{ifthen}
\newboolean{PrintVersion}
\setboolean{PrintVersion}{false}
\usepackage{mathtools}
\usepackage{comment}
\usepackage[utf8]{inputenc}
\usepackage[T1]{fontenc}
\usepackage{array}
\usepackage{tabu}
\usepackage{colortbl}
\usepackage{mathrsfs}
\usepackage{bm}

\newcommand{\mc}[1]{\mathcal{#1}}

\newcommand{\norm}[2]{\left\| #1 \right\|_{#2} }
\newcommand{\normm}[1]{\left\| #1 \right\| }
\newcommand{\inn}[2]{\left\langle #1, #2 \right\rangle}

\newcommand{\R}{\mathbb{R}}

\newcommand{{\state}}{{\mathbb{X}}}  
\newcommand{{\cs}}{{\mathbb{U}}}
\newcommand{{\as}}{{\mathbb{K}}}
\newcommand{{\Yb}}{{\mathbb{Y}}}
\newcommand{\xb}{{\bm{x}}}
\newcommand{\wb}{{\bm{w}}}
\newcommand{\vb}{{\bm{v}}}

\newcommand{\zb}{{\bm{z}}}

\newcommand{\pb}{{\bm{p}}}

\newcommand{\ub}{{\bm{u}}}
\newcommand{\rb}{{\bm{r}}}
\newcommand{\A}{{\mathcal{A}}}

\newcommand{\F}{{\mathcal{F}}}
\newcommand{\B}{{\mathcal{B}}}
\newcommand{\Bu}{ {B_{L^2 ( 0,\tau; \cs)} (R_1)}} 
\newcommand{\G}{{\mathcal{G}}}

\newcommand{\W}{{\mathbb{W}}}
\newcommand{\V}{{\mathbb{V}}}

\newcommand{\chgs}[1]{{#1}} 

\begin{document}
\begin{frontmatter}

\title{Shape Optimization of Actuators over Banach Spaces for Nonlinear Systems} 

\author[a1]{M. Sajjad Edalatzadeh}
\author[a2]{Dante Kalise}
\author[a3]{Kirsten A. Morris}
\author[a4]{Kevin Sturm}

\address[a1]{Department of Applied Mathematics, University of Waterloo, Waterloo, ON N2L 3G1, Canada (email: msedalat@uwaterloo.ca).}%
\address[a2]{School of Mathematical Sciences, University of Nottingham, University Park
Nottingham, NG7 2RD, United Kingdom (email: dante.kalise@nottingham.ac.uk).}%
\address[a3]{Department of Applied Mathematics, University of Waterloo, Waterloo, ON N2L 3G1, Canada (email: kmorris@uwaterloo.ca).}%
\address[a4]{Institute for Analysis and Scientific Computing, Technical University of Vienna, Vienna 1040, Austria (email: kevin.sturm@asc.tuwien.ac.at).}%

\begin{abstract}               
In this paper, optimal actuator shape for nonlinear parabolic systems is discussed. The system under study is an abstract differential equation with a locally Lipschitz nonlinear part. A quadratic cost on the state and input of the system is considered. The existence of an optimal actuator shape has been established in the literature. This paper focuses on driving the optimality conditions for actuator shapes belonging to a Banach space. The application of the theory to the optimal actuator shape design for railway track model is considered.
\end{abstract}


\end{frontmatter}

\section{Introduction}
Actuator shape is an important design variable for feedback synthesis in control of distributed parameter systems. Optimizing actuator shape can improve performance of the controller and significantly reduces the cost of control. Numerical simulations in \cite{kalise2017optimal} show significant improvement in the cost and performance of the control.

Optimal shape of actuators has only been studied in few works.  In (\cite{PTZ2013}), the optimal shape and position of an actuator for the wave equation in one spatial dimension are discussed. An actuator is placed on a subset $\omega\in [0,\pi]$ with a constant Lebesgue measure $L\pi$ for some $L\in (0,1)$. The optimal actuator minimizes the norm of a Hilbert Uniqueness Method (HUM)-based control; such control steers the system from a given initial state to zero state in finite time. In (\cite{privat2017actuator}), optimal actuator shape and position for  linear parabolic systems are discussed. This paper adopts the same approach as in (\cite{PTZ2013}) but with initial conditions that have randomized Fourier coefficients. The cost is defined as the average of the norm of HUM-based controls. In \cite{kalise2017optimal}, optimal actuator design for linear diffusion equations has been discussed. A quadratic cost function is considered, and shape and topological derivative of this function are derived. Optimal sensor design problems are in many ways similar to the optimal actuator design problems. In (\cite{privat2015optimal}), optimal sensor shape design has been studied where the observability is maximized over all admissible sensor shapes.Controllability-based approaches  to actuator design were used in (\cite{munch2011optimal,munch2009optimal,munch2013numerical}). 
Numerical techniques to optimize the actuator design concurrently with  the controller  are mostly limited to linear quadratic regulator problems and location of actuators, see for  example \cite{allaire2010long,kumar1978optimal-a,kubrusly1985sensors,darivandi2013algorithm}.  An $H_\infty$-approach was used in \cite{kasinathan2013h}.

The previous studies have only discussed optimal actuator shape for linear systems. Optimal actuator design problems for nonlinear distributed parameter systems has also been studied. In (\cite{edalatzadehSICON}), it is shown that under certain conditions on the nonlinearity and the cost function, an optimal input and actuator design exist, and optimality equations are derived. Results  are applied to the nonlinear railway track model as well as to the semi-linear wave model in two spatial dimensions. The existence of an optimal shape in a Banach space for nonlinear systems has been discussed in (\cite{edalatzadehTAC}). The optimality conditions in (\cite{edalatzadehTAC}) are derived  for  admissible actuator shapes in Hilbert spaces. The actuator shape space in this paper is an arbitrary Banach space.  Optimality conditions for actuator shapes over a subset of a Banach space are obtained  for nonlinear parabolic systems. A quadratic cost function on the state and input is considered to be minimized. The theory can be applied to various models. Some applications are nonlinear diffusion equation, Kuramoto-Sivashinsky equation, and nonlinear beam models (\cite{edalatzadeh2016boundary,edalatzadeh2019stability,edalatzadehSICON,edalatzadehTAC}). In this paper, the application of the theory to the optimal actuator shape design for railway track model is considered.
\section{Notation and Definitions}
Let $\state$ be a Hilbert space. The notation $\state_1\hookrightarrow \state_2$ means that the space $\state_1$ is  densely and continuously embedded in $\state_2$. Let $I$ be a set on the real line, and $m$ be a non-negative number. The Banach space $H^m(I;\state)$ is the space of all strongly measurable functions $\xb:I\to \state$ for which $\norm{\xb(t)}{\state}$ is in $H^m(I,\R)$. 
For simplicity of notation, when $I$ is an interval, the corresponding space will be indicated  without the braces; for example $L^2([0,\tau];\state)$  will be indicated by  $L^2(0,\tau;\state) . $

The Banach space $\W(0,\tau)$ is the set of all $ \xb (\cdot )  \in H^1(0,\tau;\state)\cap L^2(0,\tau;D(\A))$ with norm \cite[Section II.2]{bensoussan2015book}
\begin{equation}\notag
\norm{\xb}{\W(0,\tau)}=\norm{\dot{\xb}}{L^2(0,\tau;\state)}+\norm{\A \xb}{L^2(0,\tau;\state)}. 
\end{equation}
When there is no ambiguity, the norm on $\state$ will not be explicitly indicated. 

For every $p\in [1,\infty]$ and $\alpha\in (0,1)$, the interpolation space $D_{\A}(\alpha,p)$ is defined as the set of  all $\xb_0 \in \ss$ such that the function 
\begin{equation}
 t \mapsto v(t)\coloneqq\normm{t^{1-\alpha-1/p}\A e^{t\A}\xb_0}
\end{equation}
belongs to $L^p(0,1)$ \cite[Section 2.2.1]{lunardi2012analytic}. The norm on this space is 
$$\norm{\xb_0 }{D_{\A}(\alpha,p)}=\normm{\xb_0 }+\norm{v}{L^p(0,1)}.$$
 
\section{Optimal Actuator Design}
Let $\xb(t)$ and $\ub(t)$ be the state and input taking values in Hilbert spaces $\state$ and $\cs$, respectively. Also, let $\rb$ denote the actuator design parameter that takes value in a compact set $K_{ad}$ in a Banach space $\as$. Consider the following initial value problem (IVP):
\begin{equation}\label{eq-IVP}
\begin{cases}
\dot{\xb}(t)=\mc{A}\xb(t)+\mc{F}(\xb(t))+\mc{B}(\rb)\ub(t),\quad t>0,\\
 \xb(0)=\xb_0.
\end{cases}
\end{equation}
The nonlinear operator $\F(\cdot)$ maps a Hilbert space $\V$ to $\state$. It is assumed that $D_{\A}(1/2,2)\hookrightarrow {\V}\hookrightarrow \state.$ 

The linear operator $\A$ is associated with a sesquilinear form $a:\V\times \V\to \mathbb{C}$ (see \cite[Chapter 4]{lang2012real}). Let there be positive numbers $\alpha$ and $\beta$ such that
\begin{flalign*}
|a(\xb_1,\xb_2)|&\le \alpha \norm{\xb_1}{\V}\norm{\xb_2}{\V}, &\forall&\xb_1,\xb_2\in \V,\\
\text{Re} \; a(\xb,\xb)&\ge \beta \norm{\xb}{\V}^2,  &\forall&\xb\in \V.
\end{flalign*}
The operator $\A$ has an extension to $\bar{\A}\in\mc L(\V,\V^{^*})$ described by
\begin{equation}
\inn{\bar{\A} \vb}{\wb}_{\V^{^*},\V}=a(\vb,\wb), \quad \forall \vb, \wb \in \V,
\end{equation}
where $\V^{^*}$ denotes the dual of $\V$ with respect to pivot space $\state$.

According to \cite{edalatzadehTAC}, there are positive numbers $\tau$, $R_1$, and $R_2$ such that \eqref{eq-IVP} admits a unique solution for any initial conditions $\xb_0\in B_{\V}(R_2)$ and inputs $\ub\in\Bu$ where
\begin{flalign}\label{ad sets}
B_{\V}(R_2)&=\left\{\xb_0 \in \V: \norm{\xb_0}{\V}\le R_2   \right\},\\
\Bu&=\left\{\ub\in L^2(0,\tau;\cs): \norm{\ub}{2}\le R_1   \right\}.
\end{flalign}
The mapping $\xb=\mc S(\ub,\rb;\xb_0)$ maps input $\ub\in L^2(0,\tau;\cs)$, actuator location $\rb \in \as$ and initial condition $\xb_0\in \state$ to the corresponding solution $\xb\in \W(0,\tau)$ of \eqref{eq-IVP}.

Consider the cost function
\begin{equation*}
J(\xb,\ub)=\int_0^\tau\inn{\mc Q\xb(t)}{\xb(t)}+\inn{\mc R \ub(t)}{\ub(t)}_{\cs}dt,
\end{equation*}
where $\mc Q$ is a positive semi-definite, self-adjoint bounded linear operator on $\state$, and $\mc{R}$ is a coercive, self-adjoint linear bounded operator on $\cs$. Let $U_{ad}$ be a convex and closed set contained in the interior of $\Bu$.
For a fixed initial condition $\xb_0\in B_{\V}(R_2)$, consider the following optimization problem over the admissible input set $U_{ad}$ and actuator design set $K_{ad}$
\begin{equation}
\left\{ \begin{array}{ll}
\min&J(\xb,\ub)\\
\text{s.t.}& \xb=\mc S(\ub,\rb;\xb_0),\\
&(\ub,\rb) \in U_{ad}\times K_{ad}.
\end{array} \right. \tag{P} \label{eq-optimal problem}
\end{equation}
The existence of an optimizer to this optimization problem is proven in (\cite{edalatzadehTAC}).

\chgs{\begin{defn}
The operator $\G:\state\to \Yb$ is said to be G\^ateaux differentiable at $\xb\in \state$ in the direction $\tilde{\xb}\in \state$, if the limit
\begin{equation}
\G'(\xb;\tilde{\xb})=\lim_{\epsilon\to 0}{\frac{\|\G(\xb+\epsilon\tilde{\xb})-\G(\xb)\|_{\Yb}}{\epsilon}}
\end{equation}
exits.
\end{defn}}

The optimality conditions are derived next after assuming G\^ateaux differentiability of nonlinear operators $\F(\xb)$ and $\B(\rb)$.

\begin{enumerate}
\item[A1.] \label{as-diff F} The nonlinear operator $\F(\cdot)$ is G\^ateaux differentiable, and the derivative is linear. Indicate the G\^ateaux derivative of $\F(\cdot)$ at $\xb$ in the direction $\pb$ by ${\F}_{\xb}^\prime\pb$. Furthermore, the mapping $\xb\mapsto {\F}_{\xb}^\prime$ is bounded; that is, bounded sets in ${\V}$ are mapped to bounded sets in $\mc L({\V},\state)$.
\item[A2.] \chgs{\label{as:diff B} The control operator $\mc{B}(\rb)$ is G\^ateaux differentiable with respect to $\rb$ from $K_{ad}$ to $\mc{L}(\cs,\state)$. Indicate the G\^ateaux derivative of $\mc{B}(\rb)$ at $\rb^o$ in the direction $\tilde{\rb}$ by $\mc{B}'(\rb^o;\tilde{\rb})$. Furthermore, the mapping $\tilde{\rb}\mapsto \mc{B}'(\rb^o;\tilde{\rb})$ is bounded; that is, bounded sets in $\as$ are mapped to bounded sets in $\mc L(\cs,\state)$.}
\end{enumerate}

Using these assumptions, the G\^ateaux derivative of the solution map with respect to a trajectory $\xb(t)=\mc S(\ub(t),\rb;\xb_0)$ is calculated. The resulting map is a time-varying linear IVP. Let $\bm g\in L^p(0,\tau;\state)$, consider the time-varying system
\begin{equation}\label{eq-time var}
\begin{cases}
\dot{{\bm{h}}}(t)=(\A+{\F}_{\xb(t)}^\prime) {\bm{h}}(t)+\bm g(t),\\
{\bm{h}}(0)=0.
\end{cases}
\end{equation}

\begin{lem}\label{lem-estimate}
For every $\bm g\in L^p(0,\tau;\state)$, there is a unique solution $\bm h(t)$ to \eqref{eq-time var} in $\W(0,\tau)$. Moreover, there is a positive number $c$ independent of $\bm g$ such that
\begin{equation}\label{eq-estimate}
\|{\bm{h}}\|_{\W(0,\tau)}\le c \norm{\bm g}{L^2(0,\tau;\state)}.
\end{equation}
\end{lem}
\begin{pf}
The proof follows immediately from \cite[Corollary 5.2]{dier2015}. Let $\mc P(\cdot):[0,\tau]\to \mc L (\V,\state)$ be such that $\mc P(\cdot)\xb$ is weakly measurable for all $\xb\in \V$, and there exists an integrable function $h:[0,\tau]\to [0,\infty)$ such that $\norm{\mc P(t)}{\mc L(\V,\state)}\le h(t)$ for all $t\in [0,\tau]$. Corollary 5.2 in (\cite{dier2015}) states that for every $\xb_0\in \V$ and $\bm g\in L^2(0,\tau;\state)$, there exists a unique $\xb$ in $\W(0,\tau)$ such that
\begin{equation}
\begin{cases}
\dot{\xb}(t)=(\A  +\mc P(t)) \xb(t)+\bm g(t),\\
\xb(0)=\xb_0.
\end{cases}
\end{equation}
Moreover, there exists a constant $c>0$ independent of $\xb_0$ and $\bm g(t)$ such that
\begin{equation}
\norm{\xb}{\W(0,\tau)}^2\le c \left(\norm{\bm g}{L^2(0,\tau;\state)}^2+\norm{\xb_0}{\V}^2\right).
\end{equation} 

Since $\W(0,\tau)$ is embedded in $C(0,\tau;\V)$, the state $\xb(t)$ is bounded in $\V$ for all $t\in [0,\tau]$. This together with G\^ateaux differentiablity of $\F(\cdot)$ ensures that there is a positive number $M_{\F}$ such that
\begin{equation}
\sup_{t\in [0,\tau]}\norm{{\F}_{\xb(t)}^\prime}{\mc L(\V,\state)}\le M_{\F}.
\end{equation} 
Thus, replacing the operator $\mc P(t)$ with ${\F}_{\xb(t)}^\prime$ and noting that
\begin{equation}\label{d1}
\norm{\mc P(t)}{\mc L(\V,\state)}\le M_{\F},
\end{equation}
proves the lemma.
\end{pf}

\begin{prop}\label{prop-diff}
\chgs{Let assumptions A1 and A2 hold. The solution map $\mc S(\ub(t),\rb;\xb_0)$ is G\^ateaux differentiable with respect to each $\ub(t)$ and $\rb$ in $U_{ad}\times K_{ad}$. Let $\xb(t)=\mc S(\ub(t),\rb;\xb_0)$. The G\^ateaux derivative of $\mc S(\ub(t),\rb;\xb_0)$ at $\rb$ in the direction $\tilde{\rb}$ is the mapping $\mc S'(\ub(t),\rb;\xb_0,\tilde{\rb}) :\as\to \W(0,\tau)$, $\tilde{\rb}\mapsto \zb(t)$, where $\zb(t)$ is the strict solution to
\begin{equation}
\begin{cases}
\dot{\zb}(t)=(\A+{\F}_{\xb(t)}^\prime) \zb(t)+\B'(\rb^o;\tilde{\rb})\ub(t),\\
\zb(0)=0.
\end{cases}
\end{equation}}

\end{prop}

\begin{thm}\label{thm-optimality}
Suppose assumptions A1 and A2 hold. For any initial condition $\xb_0\in\state$, let the pair $(\ub^o,\rb^o)\in U_{ad}\times K_{ad}$ be a local minimizer of the optimization problem \eqref{eq-optimal problem} with the optimal trajectory $\xb^o=\mc{S}(\ub^o;\rb^o,\xb_0)$ and let 
 $\pb^o(t)$  indicate the strict solution in $\W(0,\tau)^*$ of the final value problem
\begin{equation}\label{adj}
\dot{\pb}^o(t)=-(\mc{A}^*+{\F_{\xb^o(t)}^\prime}^*)\pb^o(t)-\mc Q \xb^o(t), \quad \pb^o(\tau)=0.
\end{equation} 
Then $\ub^o(t)=-\mc R^{-1} \B^*(\rb^o)\pb^o(t)$ and
\begin{equation}\label{optim r}
\int_0^\tau \inn{\pb^o(t)}{\B'(\rb^o;\tilde{\rb})\ub^o(t)}dt=0
\end{equation}
for all $\tilde{\rb}\in K_{ad}$.
\end{thm}
{\bf Outline of Proof:}
Theorem 11 in (\cite{edalatzadehTAC}) ensures that $\ub^o(t)=-\mc R^{-1} \B^*(\rb^o)\pb^o(t)$. To obtain \eqref{optim r}, the G\^ateaux derivative of  $$\mc G(\ub,\rb):=J(\mc S(\ub,\rb;\xb_0),\ub)$$ at $\rb^o$ in the direction $\tilde{\rb}$ is taken. After some manipulation and integration by parts, the following G\^ateaux derivative is obtained
\begin{equation}
\mc G'(\ub^o,\rb^o;\tilde{\rb}) =\int_0^\tau \inn{\pb^o(t)}{\B'(\rb^o;\tilde{\rb})\ub^o(t)}dt.
\end{equation}
The optimality condition now follows by setting the G\^ateaux derivative to zero.

\section{Railway Track Model}
Letting $[0, \tau]$ indicate  the time interval of interest, the following semilinear PDE governs the motion of a railway track $w(x,t)$ with initial deflection $w_0(x)$ and rate of deflection $v_0(x)$ on $(x,t)\in [0,1]\times [0,\tau]$ (\cite{edalatzadehSICON}):
\begin{equation}\label{track pde}\notag
\begin{cases}
\partial_{tt}w+\partial_{xx}(\partial_{xx}w+C_d\partial_{xxt}w)+\mu \partial_t w+w+\alpha w^3\\[1mm]
\qquad=b(x,r)u(t),\\[1mm]
\allowdisplaybreaks w(x,0)=w_0(x), \quad \partial_x w(x,0)=v_0(x),\\[1mm]
\allowdisplaybreaks w(0,t)=w(1,t)=0,\\[1mm]
\allowdisplaybreaks \partial_{xx} w(0,t)+C_d\partial_{xxt}w(0,t)=0,\\[1mm]
\partial_{xx} w(1,t)+C_d\partial_{xxt} w(1,t)=0, 
\end{cases}
\end{equation}
where the subscript $\partial_x$ denotes the derivative with respect to $x$; the derivative with respect to $t$ is indicated similarly. The nonlinear part of the foundation elasticity corresponds to the coefficients $\alpha$. The constant $\mu\ge 0$ is the viscous damping coefficient of the foundation, and $C_d\ge 0$ is the coefficient of Kelvin-Voigt damping in the beam.
{The track deflection is controlled by a single external force  $u(t) . $ 
The shape influence function $b(x,r)$ is a continuous function over $[0,1]$ parametrized by the parameter $r$ that describes its dependence on the actuator design. The function $b(x,  r) $  is differentiable with respect to $r $.

Choose state $\xb:=( w,v)$ where $v=\partial_t w$   and define the state space $\state:=H^2(0,1 )\cap H_0^1(0,1 )\times L^2(0,1 )$  with  norm
\begin{equation}
\| (w,v) \|^2=\int_0^{1} (\partial_{xx} w)^2+w^2+ v^2 \, dx \label{Track-eq: norm}.
\end{equation}
Define the closed self-adjoint positive operator 
\begin{flalign}
&{A}_0w:= \partial_{xxxx} w,\notag \\
&D({A}_0):=\left\lbrace w\in H^4(0,1 )| \, w(0)=w(1)=0,\right.\notag\\
&\qquad \qquad \left. \partial_{xx}w(0)=\partial_{xx}w(1)=0 \right\rbrace,
\end{flalign}
and also define
\begin{flalign}
{A}_{\scriptscriptstyle KV}(w,v)&:=\left(v,-{A}_0(w+C_dv)\right),\\
{K}(w,v)&:=(0,-(w+\mu v)),
\end{flalign}
with 
\begin{flalign}
D({A}_{\scriptscriptstyle KV}):=&\left\lbrace(w,v)\in \state| \, v\in H^2(0,1)\cap H_0^1(0,1), \right.\notag \\
& \left. w+C_dv\in D({A}_0) \right\rbrace
\end{flalign}
The state operator $\A$ is defined as
\begin{flalign}
\allowdisplaybreaks \A  :={A}_{\scriptscriptstyle KV}+{K}, \text{ with } D ( \A  ) = D({A}_{\scriptscriptstyle KV} ).
\end{flalign}
Let $\cs:=\mathbb{R}$, the input operator $\B(r):\cs\to \state$ is
\begin{equation}
\B(r)u:=(0,b(x,r)u).
\end{equation}
The nonlinear operator $\F(\cdot):\state\to \state$ is defined as
\begin{equation}
\F(w,v):=(0,-\alpha w^3).
\end{equation}
With these definitions and by setting the state $\xb(t)=(w(\cdot,t),v(\cdot,t))$ and initial condition $\xb_0=(w_0(\cdot),v_0(\cdot))$, the state space representation of the railway track model is
\begin{equation}\label{sys}
\begin{cases}
\dot{\xb}(t)=\A \xb(t)+\F(\xb(t))+\B(r)u(t),\quad t\in (0,\tau],\\
\xb(0)=\xb_0.
\end{cases}
\end{equation}
Well-posedness and stability of this model has been established (\cite{edalatzadeh2019stability}). 

The set of admissible inputs, $U_{ad}$, is a convex and closed subset of $L^2(0,\tau)$. The set of  admissible actuator designs is denoted by $K_{ad}$ and is compact in a Banach space $\as$.

The cost function is 
\begin{equation}
J(u,r;\xb_0):=\frac{1}{2}\int_0^\tau\normm{\xb(t)}^2+\gamma \normm{u(t)}^2\; dt.
\label{eq-cost}
\end{equation}
For a fixed initial condition $\xb_0\in \state$, consider the following optimization problem over the admissible input set $U_{ad}$ and actuator design set $K_{ad}$
\begin{equation}\label{optim}
\left\{ \begin{array}{ll}
\underset{r\in K_{ad}}{\min}\;\underset{u\in U_{ad}}{\min}&J(u,r;\xb_0)\\
\text{such that}& \xb(t) \text{ solves }\eqref{sys}.
\end{array} \right. \tag{P} 
\end{equation}

Let $\pb^o(t)$  indicate the strict solution of the final value problem
\begin{equation}\label{adj}
\dot{\pb}^o(t)=-({\A}^*+{\F_{\xb^o(t)}^\prime}^*)\pb^o(t)-\frac{1}{2} \xb^o(t), \quad \xb^o(\tau)=0.
\end{equation} 
Any solution $(u^o,r^o)$ to \eqref{optim} in the interior of $U_{ad}\times K_{ad}$ satisfies 
\begin{subequations}\label{opt}
\begin{flalign}
\frac{\gamma}{2} u^o(t)+\B^*(r^o)\pb^o(t)&=0,\label{opt1}\\
\int_0^\tau \inn{\pb^o(t)}{\B'(r^o;\tilde{r})u^o(t)}dt&=0 \quad \forall \tilde{r}\in K_{ad}.\label{opt2}
\end{flalign}
\end{subequations}
Write $\pb^o(t)=(f^o(x,t),g^o(x,t))$ and $b'(x,r^o;\tilde{r})$ be the G\^ateaux derivative of $b(x,r)$ at $r^o$ in the direction $\tilde{r}$. The optimality conditions can be written
\begin{flalign}
&u^o(t)=\frac{2}{\gamma}\int_0^1b(x,r^o)g^o(x,t)dx,\\
&\int_0^\tau \int_0^1 g^o(x,t)b'(x,r^o;\tilde{r})dxdt=0, \quad \forall \tilde{r}\in K_{ad}.
\end{flalign}

\section{Future Directions}
Future directions will focus on application of the theory to various PDE models. An example of a nonlinear parabolic PDE model is Kuramoto-Sivashinsky equation which models propagation of flames as well as dynamics of thin film fluids. Another example is nonlinear models of flexible beams with Kelvin-Voigt damping. Future research will develop of suitable numerical schemes for computation of optimal actuator shapes.
\bibliography{library}
\bibliographystyle{elsarticle-num}

\end{document}